\documentclass{amsart}

\usepackage{graphicx}
\usepackage{amsmath,amssymb,amsfonts}
\usepackage{mathrsfs}
\usepackage{rotating}
\usepackage{appendix}
\usepackage[numbers]{natbib}
\usepackage{enumitem}

\numberwithin{equation}{section}

\DeclareMathOperator{\tr}{Tr}

\renewcommand{\leq}{\leqslant}
\renewcommand{\geq}{\geqslant}
\renewcommand{\Re}{\mathop{\rm{Re}}\nolimits}

\newcommand{\gU}{\mathsf{U}}

\newcommand{\gO}{\mathsf{O}}
\newcommand{\Iso}{\mathsf{Iso}}
\newcommand{\gGL}{\mathsf{GL}}

\newcommand{\R}{\mathbf{R}}
\newcommand{\C}{\mathbf{C}}
\newcommand{\K}{\mathbf{K}}

\newcommand{\M}{\mathsf{M}}

\newcommand{\scalar}[2]{\langle #1 , #2\rangle}

\newcommand{\op}{\mathrm{op}}

\newcommand{\Id}{\mathrm{Id}}

\theoremstyle{plain}

\newtheorem*{theorem*}{Theorem}

\newtheorem*{proposition*}{Proposition}

\newtheorem*{lemma*}{Lemma}

\theoremstyle{remark}

\begin{document}

\title{A characterization of inner product spaces via norming vectors}
\author{Guillaume Aubrun}
\address{Universit\'e Lyon 1, CNRS, INRIA, Institut Camille Jordan, 43, boulevard du 11
novembre 1918, 69100 VILLEURBANNE, France}
\author{Mathis Cavichioli}

\subjclass{46C15}
\keywords{Characterization of inner-product spaces, norming vectors, positive John position}

\begin{abstract}
    A finite-dimensional normed space is an inner product space if and only if the set of norming vectors of any endomorphism is a linear subspace. This theorem was proved by Sain and Paul for real scalars. In this paper, we give a different proof which also extends to the case of complex scalars.
\end{abstract}

\maketitle

\section{Introduction}

The characterization of Euclidean spaces among normed spaces, or Hilbert spaces among Banach spaces, is a classical theme in functional analysis. It can be traced back to the Jordan--von Neumann theorem \cite{JvN35}, which states that a normed space is Euclidean if and only if it satisfies the parallelogram identity. This line of research has been very fruitful: for example, the monograph \cite{Amir86} compiles about 350 characterizations of inner product spaces. 

We consider vector spaces over a field $\K$ which is either $\R$ or $\C$. Recall that an \emph{inner product} on a vector space $X$ is a map $\langle \cdot , \cdot \rangle : X \times X \to \K$ that satisfies the usual axioms of conjugate symmetry, linearity in one variable and positive-definiteness. If $\|\cdot\|$ is a norm on $X$, we say that the normed space $(X,\|\cdot\|)$ is an \emph{inner product space} if there exists an inner product $\scalar{\cdot}{\cdot}$ on $X$ such that the identity $\|x\|^2 = \scalar{x}{x}$ holds for every $x \in X$. A finite-dimensional inner product space over the real field is also called a \emph{Euclidean space}.

In this note, we consider a characterization of finite-dimensional inner product spaces by a special property of their endomorphisms. Recall that the \emph{operator norm} of a linear operator $u : X \to X$, denoted~$\| u \|_{\op}$, is defined as the smallest $C \geq 0$ such that the inequality $\|u(x) \| \leq C \|x\|$ holds for every $x \in X$. We consider the set $\mathcal{N}(u)$ of \emph{norming vectors}, defined as
\[ \mathcal{N}(u) = \{ x \in X \ : \  \|u(x)\| = \| u \|_{\op} \cdot \|x\| \}. \]

The following theorem has been proved by Sain and Paul (see~\cite[Theorem 2.2]{SainPaul13}) in the real case only. Our goal is to provide a completely different and self-contained proof, which extends naturally to the complex case.

\begin{theorem*}
Let $X$ be a finite-dimensional normed space over the real or complex field. The following are equivalent.
\begin{enumerate}[nosep]
\item The space $X$ is an inner product space.
\item For every linear operator $u:X \to X$, the set $\mathcal{N}(u)$ is a linear subspace.
\end{enumerate}
\end{theorem*}

The implication (1) $\Rightarrow$ (2) can be shown as follows. If $X$ is an inner product space, we may consider the adjoint operator $u^* : X \to X$. Let $\lambda_{\max}(u^*u)$ denote the largest eigenvalue of the self-adjoint operator $u^*u$ and $E$ be the associated eigenspace. It is elementary to check then for every $x \in X$,
\[ \|ux\|^2 = \scalar{u^*ux}{x} \leq \lambda_{\max}(u^*u) \|x\|^2 \]
with equality if and only if $x \in E$. Since $\lambda_{\max}(u^*u)=\|u\|_{\op}^2$, it follows that the set $\mathcal{N}(u)$ coincides with $E$. In particular, $\mathcal{N}(u)$ is a linear subspace.

In the following sections we prove the harder implication (2) $\Rightarrow$ (1). Our arguments are inherently finite-dimensional;  we leave open the question of whether the theorem extends to infinite-dimensional normed spaces.

\section{The main proposition}

In this section we introduce our main tool. It is strongly related to the concept of \emph{positive John position} of convex bodies which is studied in \cite{AAP22}. We use a slightly different approach that allows us to cover also the complex case in a natural way. 

Consider a norm $\|\cdot\|$ on $\K^n$ and equip the algebra $\M_n(\K)$ of $n \times n$ matrices with the corresponding operator norm $\| \cdot \|_{\op}$. We denote by $\gGL_n(\K)$ the group of invertible matrices. Given $Q \in \gGL_n(\K)$, we consider the set
\[
\mathcal{C}_Q   = \{ A \in \M_n(\K) \ : \  \| AQ \|_{\op}  \leq 1\}. 
\]
Let $\M_n^+(\K)$ be the cone of positive semi-definite matrices. The set $\M_n^+(\K) \cap \mathcal{C}_Q$ is compact and therefore contains an element of maximal determinant. (This element is unique but we do not need this information.) Such an element $A$ satisfies $\| AQ \|_{\op} = 1$.

\begin{proposition*}
Let $Q \in \gGL_n(\K)$ and $A$ of maximal determinant in $\M_n^+(\K) \cap \mathcal{C}_Q$. Then the set
\[ \mathcal{N}(AQ) = \{ x \in \K^n \ : \  \|AQx\| = \|x\| \} \]
spans $\K^n$ as a vector space.
\end{proposition*}

We show in the next section how the Proposition implies our Theorem. In the real case, the Proposition follows easily from \cite[Theorem 1.2]{AAP22}. For completeness, we include a self-contained proof.

\smallskip

{\bfseries Proof of the Proposition.}
Introduce the unit ball $\mathcal{B} = \{ x \in \K^n \ : \  \|x\| \leq 1\}$. If we identify the dual space with $\K^n$, the unit ball for the dual norm is
\[ \mathcal{B}^* = \{ y \in \K^n \ : \   | \scalar{x}{y} | \leq 1 
\textnormal{ for every } x \in \mathcal{B} 
\} .\]
By the Hahn--Banach theorem, we have for $x \in \K^n$
\[ \|x\| = \max_{y \in \mathcal{B}^*}  | \scalar{x}{y} | .\]
Consider the compact set $T = \mathcal{B} \times \mathcal{B}^*$ and let $C(T)$ be the Banach space of continuous functions from $T$ to $\K$. We define a map $\alpha : \M_n(\K) \to C(T)$ as follows: if $M \in \M_n(\K)$, $x \in \mathcal{B}$ and $y \in \mathcal{B}^*$, then
\[ \alpha(M)(x,y) = \scalar{Mx}{y}. \]
The map $\alpha$ is linear and satisfies $\| \alpha(M) \| = \| M \|_{\op}$ for every $M \in \M_n(\K)$.

Denote $\mathcal{N}=\mathcal{N}(AQ)$. Assume by contradiction that the set $\mathcal{N}$ does not span~$\K^n$. Its linear image $(A^{1/2}Q)(\mathcal{N})$ does not span $\K^n$ either, and therefore there exists a nonzero orthogonal projection $P$ such that $PA^{1/2}Qx = 0$ for every $x \in \mathcal{N}$. Let $H$ be the matrix defined as $H= \lambda P- \Id$, where $\lambda$ is a positive number chosen so that $\tr H >0$.

Introduce functions $f$, $g \in C(T)$ defined as $f = \alpha(AQ)$ and $g = \alpha(A^{1/2}HA^{1/2}Q)$. Observe that $\|f\|=\| AQ \|_{\op} =1$. For $t=(x,y) \in T$, we have a chain of implications
\[ |f(t)|=1 \Rightarrow \|AQ(x)\|=1 \Rightarrow x \in \mathcal{N} \Rightarrow HA^{1/2}Qx = -A^{1/2}Qx \Rightarrow g(t)=-f(t).\] 
Consequently, the functions $f$ and $g$ satisfy the hypothesis of the following lemma, whose proof is postponed.

\begin{lemma*}
Let $T$ be a nonempty compact topological space and $f,g \in C(T)$. Assume that for every $t \in T$ such that $|f(t)|=\|f\|$, we have $\Re f(t)\overline{g(t)}  <0$. Then, for $\delta>0$ small enough, we have $\|f + \delta g\| < \|f\|$.
\end{lemma*}

The lemma implies that for $\delta >0$ small enough, we have
\[ \| A^{1/2}(\Id + \delta H)A^{1/2}Q \|_{\op} = \|f + \delta g\| \leq 1 \]
and thus $A^{1/2}(\Id + \delta H)A^{1/2} \in \mathcal{C}_Q$. 
As $\delta$ goes to zero, we have $\det(\Id + \delta H) = 1 + \delta \tr H+ o(\delta)$ and therefore $\det(\Id + \delta H)>1$ for $\delta >0$ small enough. The inequality
\[ \det(A^{1/2}(\Id + \delta H)A^{1/2} ) = \det(A)\det( \Id + \delta H) > \det(A)\]
contradicts the maximality of $A$. \hfill\qedsymbol

\medskip

It remains to prove the lemma.

\smallskip

{\bfseries Proof of the Lemma.}
Both $f$ and $g$ are nonzero (otherwise the hypothesis fails) and we may assume by rescaling that $\|f\|=\|g\|=1$.
Let $T_1$ be the nonempty closed subset of $T$ defined as $T_1 = \{t \in T \ : \  |f(t)|=1 \}$. Since the function $t \mapsto \Re f(t)\overline{g(t)}$ is continuous, it achieves its maximum on $T_1$ and therefore there exists $\varepsilon >0$ such that $\Re f(t)\overline{g(t)} \leq -\varepsilon$ for every $t \in T_1$. Denote by $T_2$ the closed subset of $T$ defined as 
$T_2 = \{t \in T \ : \  \Re f(t)\overline{g(t)} \geq -\varepsilon \}$. Since $f$ is continuous, there exists $\eta>0$ such that $|f(t)| \leq 1-\eta$ for every $t \in T_2$.
For $t \in T$ and $\delta>0$, we compute
\begin{eqnarray*} |(f+\delta g)(t)|^2 &= &|f(t)|^2 + 2 \delta \Re f(t) \overline{g(t)} + \delta^2 |g(t)|^2 \\
& \leq & 
\begin{cases}
(1-\eta)^2 + 2 \delta + \delta^2 & \textnormal{if } t \in T_2\\
1 - 2 \delta \varepsilon + \delta^2  & \textnormal{if } t \not \in T_2
\end{cases}
\end{eqnarray*}
It follows that $\| f +\delta g\|^2 \leq \max((1-\eta)^2 + 2 \delta + \delta^2,1 - 2 \delta \varepsilon + \delta^2)$ and therefore $\| f +\delta g\|<1$ for $\delta >0$ small enough. \hfill\qedsymbol

\section{Proof of the Theorem}

The implication (1) $\Longrightarrow$ (2) has been proved in the introduction. Conversely, let~$X$ be a finite-dimensional normed space satisfying condition (2) from the Theorem. Let $\Iso(X)$ be the group of isometries of $X$, defined as
\[ \Iso(X) = \{ u : X \to X \textnormal{ linear } \ : \  \|u(x)\|=\|x\| \textnormal{ for every }x \in X \}.  \]
There is an inner product on $X$ which is invariant with respect to $\Iso(X)$ (see for example \cite[p.~131, Theorem 2]{OV90}) and therefore, without loss of generality, we may assume that $X=(\R^n,\|\cdot\|)$ and that $\Iso(X)$ is a subgroup of the orthogonal group~$\gO_n$. (From now on we consider only the real case, the complex case is similar using the unitary group~$\gU_n$.)

Fix $Q \in \gO_n$ and let $A$ be an element of maximal determinant in $\M_n^+(\K) \cap \mathcal{C}_Q$. The set $\mathcal{N}(AQ)$ is a linear subspace of $\R^n$ (by hypothesis (2) from the Theorem) which spans~$\R^n$ (by the conclusion of the Proposition). It follows that $\mathcal{N}(AQ) = \R^n$ and therefore that $AQ \in \Iso(X) \subset \gO_n$. The matrix $A = (AQ)Q^{-1}$ is both orthogonal and positive semidefinite; it follows that $A$ is the identity matrix and therefore $Q \in \Iso(X)$.

The previous paragraph shows that $\Iso(X) = \gO_n$. As a consequence, for every~$x$ and $y$ in the Euclidean unit sphere $S^{n-1}$, there exists $u \in \Iso(X)$ such that $u(x)=y$. It follows that the norm $\|\cdot\|$ is constant on $S^{n-1}$, hence is a multiple of the Euclidean norm. This shows that $X$ is an inner product space.

\section*{Acknowledgement}

The first author was supported in part by ANR under the grant ESQuisses
(ANR-20-CE47-0014-01). We thank Miguel Mart\'in for pointing to us the reference~\cite{SainPaul13}.

\bibliography{references}{}
\bibliographystyle{plain}

\end{document}